\documentclass[a4paper,11pt,english]{smfart}

\usepackage{amssymb}
\usepackage{url}
\usepackage[T1]{fontenc}
\RequirePackage{calrsfs}
\DeclareSymbolFont{rsfscript}{OMS}{rsfs}{m}{b}
\DeclareSymbolFontAlphabet{\mathrsfs}{rsfscript}
\usepackage[utopia,expert]{mathdesign}
\usepackage{lscape}

\usepackage[leftbars]{changebar}

\usepackage{palatino}
\usepackage{rotating}
\usepackage{graphicx}

\usepackage{enumerate}

\usepackage{color}
\definecolor{shadecolor}{gray}{0.90}

\def\bfit{\bfseries\itshape}

\def\proofname{D\'emonstration}

\input xypic
\xyoption{all}
\xyoption{arc}

\makeindex

\def\indexname{Index des notations}
\newtheorem{theo}{Theorem}[section]
\def\thetheo{\thesection.\arabic{theo}}
\newtheorem{prop}[theo]{Proposition}

\newtheorem{lem}[theo]{Lemma}

\newtheorem{coro}[theo]{Corollary}

\newtheorem{defi}[theo]{Definition}

\renewcommand\theequation{\thetheo}
\def\equat{\refstepcounter{theo}\begin{equation}}
\def\endequat{\end{equation}}

\renewcommand\thetable{\thetheo}

\renewcommand\thesection{\arabic{section}}
\renewcommand\thesubsection{\thesection.\Alph{subsection}}

\def\contentsname{Table des mati\`eres}
\def\listfigurename{Liste des figures}
\def\listtablename{Liste des tables}
\def\appendixname{Appendice}

\def\refname{R\'ef\'erences}

    \def\CM{{\mathbb{C}}}

    \def\SM{{\mathbb{S}}}
    \def\TM{{\mathbb{T}}}

    \def\ZM{{\mathbb{Z}}}

    \def\CC{{\mathcal{C}}}
    
    \def\EC{{\mathcal{E}}}
    \def\FC{{\mathcal{F}}}

  \def\kb{{\mathbf k}}  
    \def\LC{{\mathcal{L}}}
  \def\mb{{\mathbf m}}  \def\MC{{\mathcal{M}}}

  \def\rb{{\mathbf r}}  
\def\Sb{{\mathbf S}}    \def\SC{{\mathcal{S}}}
\def\Tb{{\mathbf T}}    
    
  \def\vb{{\mathbf v}}  
    
\def\Xb{{\mathbf X}}    
    
\def\Zb{{\mathbf Z}}

    \def\ZCB{{\boldsymbol{\mathcal{Z}}}}

\def\a{\alpha}
\def\b{\beta}

\def\d{\delta}
\def\D{\Delta}
\def\e{\varepsilon}
\def\ph{\varphi}
\def\l{\lambda}
\def\L{\Lambda}

\def\o{\omega}
\def\O{\Omega}

\def\th{\theta}

\def\t{\tau}

\def\z{\zeta}

\def\gamb{{\boldsymbol{\gamma}}}

        \def\pht{{\tilde{\varphi}}}

\def\mub{{\boldsymbol{\mu}}}

\DeclareMathOperator{\diag}{{\mathrm{diag}}}

\DeclareMathOperator{\End}{{\mathrm{End}}}
\DeclareMathOperator{\Ext}{{\mathrm{Ext}}}

\DeclareMathOperator{\Id}{{\mathrm{Id}}}

\DeclareMathOperator{\im}{{\mathrm{Im}}}

\DeclareMathOperator{\Irr}{{\mathrm{Irr}}}
\DeclareMathOperator{\Ker}{{\mathrm{Ker}}}

\DeclareMathOperator{\Rad}{{\mathrm{Rad}}}
\DeclareMathOperator{\Res}{{\mathrm{Res}}}

\DeclareMathOperator{\Tr}{{\mathrm{Tr}}}

\def\to{\rightarrow}
\def\longto{\longrightarrow}
\def\injto{\hookrightarrow}

\def\fonction#1#2#3#4#5{\begin{array}{rccc}
{#1} : & {#2} & \longto & {#3}  \\
& {#4} & \longmapsto & {#5} 
\end{array}}

\def\fonctio#1#2#3#4{\begin{array}{ccc}
{#1} & \longto & {#2} \\
{#3} & \longmapsto & {#4} 
\end{array}}

\def\DS{\displaystyle}

\def\finl{~$\blacksquare$}

\def\lexp#1#2{\kern\scriptspace\vphantom{#2}^{#1}\kern-\scriptspace#2}
\def\le{\hspace{0.1em}\mathop{\leqslant}\nolimits\hspace{0.1em}}
\def\ge{\hspace{0.1em}\mathop{\geqslant}\nolimits\hspace{0.1em}}

\mathchardef\inferieur="321E
\mathchardef\superieur="321F

\def\eqna{\begin{eqnarray*}}
\def\endeqna{\end{eqnarray*}}

\def\itemth#1{\item[${\mathrm{(#1)}}$]}

\catcode`\@=11
\long\def\@car#1#2\@nil{#1}
\long\def\@first#1#2{#1}
\long\def\@second#1#2{#2}
\long\def\ifempty#1{\expandafter\ifx\@car#1@\@nil @\@empty
  \expandafter\@first\else\expandafter\@second\fi}
\catcode`\@=12

\def\boitegrise#1#2{\begin{centerline}{\fcolorbox{black}{shadecolor}{~
    \begin{minipage}[t]{#2}{\vphantom{~}#1\vphantom{$A_{\DS{A_A}}$}}
            \end{minipage}~}}\end{centerline}\medskip}

\def\surto{\twoheadrightarrow}
\def\module{{\text{-}}{\mathrm{mod}}}

\theoremstyle{remark}

\theoremstyle{plain}

\def\BIL{LR}
\def\GAUCHE{L}
\def\CAR{CAR}
\def\FAM{FAM}

\def\xyinj{\ar@{^{(}->}}
\def\xysur{\ar@{->>}}

\def\isomorphisme#1{{\boldsymbol{[}}\hskip0.5mm #1\hskip0.5mm {\boldsymbol{]}}}
\def\res{{\mathrm{res}}}

\bigskip\def\unb{{\boldsymbol{1}}}

\def\trace{{\boldsymbol{t}}_{\!\SSS{\HC}}}

\DeclareMathOperator{\ev}{{\mathrm{ev}}}
\DeclareMathOperator{\coev}{{\mathrm{coev}}}

\DeclareMathOperator{\groth}{{\mathrm{Gr}}}

\makeatletter
\def\hlinewd#1{%
\noalign{\ifnum0=`}\fi\hrule \@height #1 %
\futurelet\reserved@a\@xhline}
\makeatother

\newlength\epaisLigne

\usepackage{multirow}

\newcommand{\longiso}{\stackrel{\sim}{\longrightarrow}}

\def\modules{\operatorname{\!-mod}\nolimits}

\def\stable{\operatorname{\!-stab}\nolimits}
\def\STABLE{\operatorname{\!-stab_B}\nolimits}
\def\projectifs{\operatorname{\!-proj}\nolimits}
\def\PROJECTIFS{\operatorname{\!-proj_B}\nolimits}

\def\sta{{\mathrm{st}}}

\makeatletter
\def\hlinewd#1{%
\noalign{\ifnum0=`}\fi\hrule \@height #1 %
\futurelet\reserved@a\@xhline}
\makeatother

\usepackage{multirow}

\usepackage{lscape}

\usepackage{pstricks}

\begin{document}

%\baselineskip=16pt
%\large\baselineskip=20pt
%\Large\baselineskip=24pt

\title{An asympotic cell category for cyclic groups}

\author{{\sc C\'edric Bonnaf\'e}}
\address{
Institut Montpelli\'erain Alexander Grothendieck (CNRS: UMR 5149), 
Universit\'e Montpellier 2,
Case Courrier 051,
Place Eug\`ene Bataillon,
34095 MONTPELLIER Cedex,
FRANCE} 

\makeatletter
\email{cedric.bonnafe@umontpellier.fr}
\makeatother

\author{{\sc Rapha\"el Rouquier}}

\address{UCLA Mathematics Department
Los Angeles, CA 90095-1555, 
USA}
\email{rouquier@math.ucla.edu}

\date{\today}

\thanks{The first author is partly supported by the ANR 
(Project No ANR-16-CE40-0010-01 GeRepMod). \\
The second author was partly 
 supported by the NSF (grant DMS-1161999 and DMS-1702305) and by a
grant from the Simons Foundation (\#376202)}
\maketitle
\pagestyle{myheadings}
\markboth{\sc C. Bonnaf\'e \& R. Rouquier}{\sc An asympotic cell category for cyclic groups}

\begin{center}
To Michel Brou\'e, for his spetsial intuitions.
\end{center}

\vskip1cm

\begin{abstract}
In his theory of unipotent characters of finite groups of Lie type, Lusztig
constructed modular categories from two-sided cells in Weyl groups. Brou\'e,
Malle and Michel have extended parts of Lusztig's theory to complex reflection groups.
This includes generalizations of the corresponding fusion algebras, although the
presence of negative structure constants prevents them from arising from modular
categories. We give here the first construction of braided pivotal monoidal categories
associated with non-real reflection groups (later reinterpreted by Lacabanne as
super modular categories). They are associated with
cyclic groups, and their fusion algebras are those constructed by Malle.
\end{abstract}

\bigskip

Brou\'e, Malle and Michel \cite{bmm}
have constructed a combinatorial version of Lusztig's theory
of unipotent characters of finite groups of Lie type for certain complex reflection
groups ("spetsial groups").
The case of spetsial imprimitive complex reflection groups has been
considered by Malle in~\cite{malle}: Malle defines a combinatorial set
which generalizes the one defined by Lusztig to parametrize 
{\it unipotent characters} of the associated 
finite reductive group when $W$ is a Weyl group. Malle generalizes also the partition
of this set into Lusztig {\it families}. To each family, he associates a
{$\ZM$-fusion datum}: 
a $\ZM$-fusion datum is a structure similar to a 
usual fusion datum (which we will call a {\it $\ZM_+$-fusion datum})
except that the structure constants of the associated 
fusion ring might be negative. 

It is a classical problem to find a 
tensor category with suitable extra-structu- -res (pivot, twist) corresponding
to a given $\ZM_+$-fusion datum.
The aim of this paper is to provide an {\it ad hoc} categorification of the 
$\ZM$-fusion datum associated with the non-trivial family of 
the cyclic complex reflection group of order $d$: it is provided 
by a quotient category of the representation category of the
Drinfeld double of the {\it Taft algebra} of dimension $d^2$ 
(the Taft algebra is the positive part of a Borel subalgebra of quantum 
$\mathfrak{sl}_2$ at a $d$-th root of unity). 

The constructions of this paper have recently been extended by 
Lacabanne \cite{abel 1},~\cite{abel 2} to some families of the complex reflection 
groups $G(d,1,n)$.
His construction involves super-categories instead 
of triangulated categories, as suggested by Etingof. In the particular case studied 
in our paper, Lacabanne's construction gives a reinterpretation 
as well as some clarifications of our results.

\bigskip
In the first section, we recall some basic properties of the Taft algebra and
its Drinfeld double. The second section is devoted to recalling some of
the structure of its category of representations: simple modules, blocks and
structure of projective modules. We summarize some elementary facts on the tensor
structure in the third section: invertible objects and tensor product of simple
objects by the defining two-dimensional representation.
The fourth section provides generators and characters of
the Grothendieck group of $D(B)$, a commutative ring.

Our original work starts in the fifth section, with the study of the
stable module category of $D(B)$ and a further quotient and the
determination of their Grothendieck rings. In the sixth section, we define and
study a pivotal structure on $D(B)\modules$ and determine characters of
its Grothendieck group associated to left or right traces. This gives
rise to positive and negative $S$ and $T$-matrices. We proceed similarly for
the small quotient triangulated category. These are our fusion data.

We recall in the final section the construction of Malle's fusion datum 
associated with cyclic groups and we check that it coincides with the
fusion data we defined in the previous section.

In Appendix~\ref{sec:preliminaires}, we provide a description of $S$-matrices in
the setting of pivotal tensor categories.

\bigskip

\noindent{\sc Acknowledgements - } We wish to thank warmly Abel Lacabanne for 
the fruitful discussions we had with him and 
his very careful reading of a preliminary version of this work.

\bigskip

% \subsection{Twist} 
% In this section, we also assume that $\CC$ is endowed with a {\it twist} $\th$ 
% (see???), that is, a functorial family of isomorphisms $\th_X : X \longiso X$ 
% such that 
% $$\th_{X \otimes Y} = (\th_X \otimes \th_Y) c_{Y,X} c_{X,Y}$$
% for all objects $X$, $Y$ in $\CC$ (see~\cite[Definition~8.10.1]{egno}). Note that 
% we do not assume that $\th$ is a {\it ribbon structure} (as well as we do 
% not assume that $\CC$ is spherical). If $X$ is endominimal, we view 
% $\th_X \in \End_\CC(X)=\CM$ as a scalar. The next result is ??????????:
% 
% \bigskip
% 
% \begin{prop}\label{prop:theta}
% Let $X$ and $Y$ be two endominimal objects in $\CC$. Then
% $$\dim_-(X) s_{X,Y}^+ = \th_X^{-1}\th_Y^{-1} \sum_{Z \in \Irr(\CC)} [X \otimes Y : Z] 
% \th_Z \dim_-(Z).$$
% \end{prop}
% 
% \bigskip

\section{The Drinfeld double of the Taft algebra}

\medskip

From now on, $\otimes$ will denote the tensor product $\otimes_\CM$. 
We fix a natural number $d \ge 2$ as well as a primitive $d$-th root 
of unity $\z \in \CM^\times$. We denote by $\mub_d=\langle \z \rangle$ 
the group of $d$-th roots of unity. 
%Given $i \in \ZM/d\ZM$ and $\a$  an element 

%of a group such that $\a^d=1$, we will denote by $\a^i$ the element $\a^{i'}$, where 

%$i'$ is any representative of $i$.

Given $n \ge 1$ a natural number and $\xi \in \CM$, we set
$$(n)_\xi=1+\xi+\cdots + \xi^{n-1}$$
$$(n)!_\xi=\prod_{i=1}^n (i)_\xi.\leqno{\text{and}}$$
We also set $(0)!_\xi=1$.

\bigskip

\subsection{The Taft algebra}
We denote by $B$ the $\CM$-algebra defined by the following presentation:
\begin{itemize}
\item[$\bullet$] Generators: $K$, $E$.

\item[$\bullet$] Relations:
$\begin{cases}
K^d=1,\\
E^d=0,\\
KE=\z EK.
\end{cases}$
\end{itemize}
It follows from~\cite[Proposition~IX.6.1]{kassel} that:
\begin{itemize}
\itemth{\D} There exists a unique morphism of algebras $\D : B \to B \otimes B$ 
such that
$$\D(K)=K \otimes K\qquad\text{and}\qquad \D(E)=(1 \otimes E) + (E \otimes K).$$

\itemth{\e} There exists a unique morphism of algebras $\e : B \to \CM$ such that 
$$\e(K)=1\qquad\text{and}\qquad \e(E)=0.$$

\itemth{S} There exists a unique anti-automorphism $S$ of $B$ such that
$$S(K)=K^{-1}\qquad\text{and} \qquad S(E)=-EK^{-1}.$$
\end{itemize}
With $\D$ as a coproduct, $\e$ as a counit and $S$ as an (invertible) antipode, 
$B$ becomes a Hopf algebra, called the {\it Taft algebra}~\cite[Example~5.5.6]{egno}. 
It is easily checked that 
\equat\label{eq:base-B}
B=\bigoplus_{i,j =0}^{d-1} \CM K^i E^j = \bigoplus_{i,j=0}^{d-1} \CM E^iK^j.
\endequat

\bigskip

\subsection{Dual algebra} 
Let $K^*$ and $E^*$ denote the elements of $B^*$ 
such that 
$$K^*(E^iK^j)=\d_{i,0} \z^j\qquad\text{and}\qquad E^*(E^iK^j)=\d_{i,1}.$$
Recall that $B^*$ is naturally a Hopf algebra~\cite[Proposition~III.3.3]{kassel} 
and it follows from~\cite[Lemma~IX.6.3]{kassel} that
\equat\label{eq:mult-dual}
(E^{*i} K^{*j})(E^{i'}K^{j'})=\d_{i,i'} (i)!_\z \z^{j(i+j')}.
\endequat
We deduce easily that $(E^{*i} K^{*j})_{0 \le i,j \le d-1}$ is a $\CM$-basis 
of $B^*$. 

\def\cop{{\mathrm{cop}}}

We will give explicit formulas for the coproduct, the counit and the antipode 
in the next subsection.  We will in fact use the Hopf algebra $(B^*)^\cop$, 
which is the Hopf algebra whose underlying space is $B^*$, 
whose product is the same as in $B^*$ and whose coproduct 
is opposite to the one in $B^*$. 

\bigskip

\subsection{Drinfeld double} 
We denote by $D(B)$ the {\it Drinfeld quantum double} of $B$, as defined 
for instance in~\cite[Definition~IX.4.1]{kassel} or~\cite[Definition~7.14.1]{egno}. 
Recall that $D(B)$ contains $B$ and $(B^*)^\cop$ as Hopf subalgebras and 
that the multiplication induces an isomorphism of vector spaces 
$(B^*)^\cop \otimes B \longiso D(B)$. 
A presentation of $D(B)$, with generators $E$, $E^*$, 
$K$, $K^*$ is given for instance in~\cite[Proposition~IX.6.4]{kassel}. 
We shall slightly modify it by setting
$$z=K^{*-1}K\qquad\text{and}\qquad F=\z E^* K^{*-1}.$$
Then~\cite[Proposition~IX.6.4]{kassel} can we rewritten as follows:

\bigskip

\begin{prop}\label{prop:presentation}
The $\CM$-algebra $D(B)$ admits the following presentation:
\begin{itemize}
\item[$\bullet$] Generators: $E$, $F$, $K$, $z$;

\item[$\bullet$] Relations: 
$\begin{cases}
K^d=z^d=1,\\
E^d=F^d=0,\\
[z,E]=[z,F]=[z,K]=0,\\
KE=\z EK,\\
KF=\z^{-1} FK,\\
[E,F]=K-z K^{-1}.
\end{cases}$
\end{itemize}
\end{prop}

\bigskip

The next corollary follows from an easy induction argument:

\bigskip

\begin{coro}\label{coro:fei}
If $i \ge 1$, then 
$$[E,F^i]=(i)_\z F^{i-1}(\z^{1-i}K-zK^{-1})$$
$$[F,E^i]=(i)_\z E^{i-1}(\z^{1-i}zK^{-1}-K).\leqno{\text{and}}$$
\end{coro}

\bigskip

The algebra $D(B)$ is endowed with a structure of Hopf algebra, where the 
comultiplication, the counit and the antipode are still denoted by $\D$, $\e$ 
and $S$ respectively (as they extend the corresponding objects for $B$). We 
have~\cite[Proposition~IX.6.2]{kassel}:
\equat\label{eq:delta-s}
\begin{cases}
\D(K)=K \otimes K,\\
\D(z)=z \otimes z,\\
\D(E)=(1 \otimes E) + (E \otimes K),\\
\D(F)=(F \otimes 1) + (zK^{-1} \otimes F),\\
\end{cases}
\qquad
\begin{cases}
S(K)=K^{-1},\\
S(z)=z^{-1},\\
S(E)=-EK^{-1},\\
S(F)=-\z^{-1}FKz^{-1},\\
\end{cases}
\endequat
\equat\label{eq:counit}
\e(K)=\e(z)=1\qquad\text{and}\qquad\e(E)=\e(F)=0.
\endequat

\bigskip

\subsection{Morphisms to $\CM$}\label{sub:morphismes}
Given $\xi \in \mub_d$, we denote by $\e_\xi : D(B) \to \CM$ the unique  
morphism of algebras such that 
$$\e_\xi(K)=\xi,\quad\e_\xi(z)=\xi^2\quad\text{and}\quad\e_\xi(E)=\e_\xi(F)=0.$$
It is easily checked that the $\e_\xi$'s are the only morphisms of algebras 
$D(B) \to \CM$. Note that $\e_1=\e$ is the counit. 

\bigskip

\subsection{Group-like elements}
It follows from~(\ref{eq:delta-s}) that $K$ and $z$ are group-like, so that 
$K^iz^j$ is group-like for all $i$, $j \in \ZM$. The converse also holds
(and is certainly already well-known).

\bigskip

\begin{lem}\label{lem:group-like}
If $g \in D(B)$ is group-like, then there exist $i$, $j \in \ZM$ 
such that $g=K^iz^j$.
\end{lem}

\bigskip

\begin{proof}
Let $g \in D(B)$ be a group-like element. Let us write 
$$g=\sum_{i,j,k,l=0}^{d-1} \a_{i,j,k,l} K^i z^j E^k F^l.$$
We denote by $(k_0,l_0)$ the biggest pair (for the lexicographic order) 
such that there exist $i$, $j \in \{0,1,\dots,d-1\}$ such that $\a_{i,j,k_0,l_0} \neq 0$. 
The coefficient of $K^iz^jE^{k_0} F^{l_0} \otimes K^iz^jE^{k_0} F^{l_0}$ 
in $g \otimes g$ is equal to $\a_{i,j,k_0,l_0}^2$, so it is different from $0$. 

But, if we compute the coefficient of $K^iz^jE^{k_0} F^{l_0} \otimes K^iz^jE^{k_0} F^{l_0}$ in 
$$g \otimes g = \D(g)=\sum_{i,j,k,l=0}^{d-1} \a_{i,j,k,l} \D(K)^i \D(z)^j \D(E)^k \D(F)^l$$
using the formulas~(\ref{eq:delta-s}), we see that it is equal to $0$ 
if $(k_0,l_0) \neq (0,0)$. Therefore $(k_0,l_0) = (0,0)$, and so 
$g$ belongs to the linear span of the family $(K^iz^j)_{i,j \in \ZM}$. 
Now the result follows from the linear independence of group-like elements.
\end{proof}

\bigskip

\bigskip

\subsection{Braiding}\label{sub:braiding}
For $0 \le i,j \le d-1$, we set
$$\b_{i,j}=\frac{E^{*i}}{d\cdot (i)!_\z} \sum_{k=0}^{d-1} \z^{-k(i+j)} K^{*k}.$$
It follows from~(\ref{eq:mult-dual}) that $(\b_{i,j})_{0 \le i,j \le d-1}$ is a 
dual basis to $(E^iK^j)_{0\le i,j \le d-1}$. 
We set now
$$R=\sum_{i,j=0}^{d-1} E^iK^j \otimes \b_{i,j} \in D(B) \otimes D(B).$$
Note that $R$ is a universal $R$-matrix for $D(B)$ and it endows $D(B)$ with a structure of 
braided Hopf algebra~\cite[Theorem~IX.4.4]{kassel}). 
Using our generators $E$, $F$, $K$, $z$, we have:
\equat\label{eq:r-matrice}
R=\frac{1}{d}\sum_{i,j,k=0}^{d-1} 
\frac{\z^{(i-k)(i+j)-i(i+1)/2}}{(i)!_\z} E^iK^j \otimes z^{-k}F^iK^k.
\endequat

\bigskip

\subsection{Twist}\label{sub:theta}
Let us define
$$\fonction{\t}{D(B) \otimes D(B)}{D(B) \otimes D(B)}{a \otimes b}{b \otimes a.}$$
Following~\cite[\S{VIII.4}]{kassel}, we set 
$$u=\sum_{i,j=0}^{d-1} S(\b_{ij}) E^iK^j \in D(B).$$
Recall that $u$ is called the {\it Drinfeld element} of $D(B)$. 
It satisfies several properties (see for instance~\cite[Proposition~VIII.4.5]{kassel}). 
For instance, $u$ is invertible and 
we will recall only three equalities:
\equat\label{eq:drinfeld}
\e(u)=1,\quad \D(u)=(\t(R)R)^{-1}(u \otimes u)\quad\text{and}\quad
S^2(b)=ubu^{-1}
\endequat
for all $b \in D(B)$. A straightforward computation shows that 
\equat\label{eq:s-carre}
S^2(b)=KbK^{-1}
\endequat
for all $b \in D(B)$. We now set
$$\th=K^{-1}u.$$
The following proposition is a consequence of~(\ref{eq:drinfeld}) and~(\ref{eq:s-carre}).

\bigskip

\begin{prop}\label{prop:twist}
The element $\th$ is central and invertible in $D(B)$ and satisfies
$$\e(\th)=1\qquad\text{and}\qquad \D(\th)=(\t(R)R)^{-1}(\th \otimes \th).$$
\end{prop}

\bigskip

Let us give a formula for $\th$:

\equat\label{eq:theta}
\th=\frac{1}{d}\sum_{i,j,k =0}^{d-1} (-1)^i 
\frac{\z^{(i-k)(i+j)-i}}{(i)!_\z} z^{k-i}F^iE^iK^{i+j-k-1}.
\endequat

\bigskip

\begin{coro}\label{coro:s-theta}
We have $S(\th)=z\th$.
\end{coro}

\bigskip

\begin{proof}
Let $g=S(\th)\th^{-1}$. 
Since $\D \circ S = \t \circ (S \otimes S) \circ \D$ and $(S \otimes S)(R)=R$ 
(see for instance~\cite[Theorems~III.3.4~and~VIII.2.4]{kassel}, it follows from 
Proposition~\ref{prop:twist} that $g$ is central and group-like. 
Hence, by Lemma~\ref{lem:group-like}, there exists $l \in \ZM$ such that 
$S(\th)=\th z^l$. 
So, by~(\ref{eq:theta}), we have
$$S(\th)E^{d-1}=\th z^l E^{d-1}=
\frac{1}{d} \sum_{j,k \in \ZM/d\ZM}^{d-1} \z^{-jk} z^{k+l} K^{j-k-1} E^{d-1}.\leqno{(\sharp)}$$
Let us now compute $S(\th) E^{d-1}$ by using directly~(\ref{eq:theta}). We get
\eqna
S(\th)E^{d-1}=E^{d-1}S(\th)&=&\DS{\frac{1}{d} \sum_{j,k \in \ZM/d\ZM} 
\z^{-jk} E^{d-1} z^{-k} K^{1+k-j}}\\
&=&\DS{\frac{1}{d} 
\sum_{j,k \in \ZM/d\ZM} \z^{-jk} \z^{1+k-j} z^{-k} K^{1+k-j} E^{d-1}}\\
&=&\DS{\frac{1}{d} 
\sum_{j,k \in \ZM/d\ZM} \z^{(1-j)(1+k)} z^{-k} K^{1+k-j} E^{d-1}}.
\endeqna
So, if we set $j'=1-j$ and $k'=-1-k$, we get
$$S(\th)E^{d-1} = \frac{1}{d} \sum_{j',k' \in \ZM/d\ZM} \z^{-j'k'}z^{k'+1}K^{j'-k'-1} E^{d-1}.$$
Comparing with $(\sharp)$, we get that $z^l=z$.
\end{proof}

\bigskip

\section{$D(B)$-modules} 

\medskip

Most of the result of this section are due to Chen~\cite{chen 1} 
or~Erdmann, Green, Snashall and Taillefer~\cite{egst1},~\cite{egst2}.
By a $D(B)$-module, we mean a finite dimensional left $D(B)$-module. 
We denote by $D(B)\modules$ the category of (finite dimensional left) $D(B)$-modules. 
% If $l \ge 1$ and 
% $1 \le i,j \le l$, we denote by $\EM_{i,j}^{(l)}$ the elementary matrix 
% whose non-zero entry is on the $i$-th row and the $j$-th column. 
Given $\a_1$,\dots, $\a_{l-1} \in \CM$, we set
$$J_l^+(\a_1,\dots,\a_{l-1})=
% \a_1 \EM_{1,2}{(l)} + \a_2 \EM_{2,3}^{(l)} 
% + \cdots + \a_{l-1} \EM_{l-1,l}^{(l)}=
\begin{pmatrix}
0 & \a_1 & 0 & \cdots & 0 \\
\vdots & \ddots & \ddots & \ddots &\vdots \\
\vdots & & \ddots & \ddots & 0 \\
0 &  &  & 0 & \a_{l-1} \\
0 & \cdots & \cdots & \cdots & 0 \\
\end{pmatrix}$$
$$J_l^-(\a_1,\dots,\a_{l-1})=\lexp{t}{J_l^+(\a_1,\dots,\a_{l-1})}.\leqno{\mathrm{and}}$$
Given $M$  a $D(B)$-module and $b \in D(B)$, we denote by $b|_M$ the endomorphism of 
$M$ induced by $b$. For instance, $E|_M$ and $F|_M$ are nilpotent and 
$K|_M$ and $z|_M$ are semisimple.

\bigskip

\subsection{Simple modules}
Given $1 \le l \le d$ and $p \in \ZM/d\ZM$, we denote by $M_{l,p}$ the 
$D(B)$-module with $\CM$-basis $\MC^{(l,p)}=(e_i^{(l,p)})_{1 \le i \le l}$ where
the action of $z$, $K$, $E$ and $F$ in the basis $\MC^{(l,p)}$ 
its given by the following matrices:
$$z|_{M_{l,p}}=\z^{2p+l-1} \Id_{M_{l,p}},$$
$$K|_{M_{l,p}}=\z^p \diag(\z^{l-1},\z^{l-2},\dots,\z,1),$$
$$E|_{M_{l,p}}=\z^p J_l^+((1)_\z(\z^{l-1}-1),(2)_\z(\z^{l-2}-1),
\dots,(l-1)_\z(\z-1)),$$
$$F|_{M_{l,p}}=J_l^-(1,\dots,1).$$
It is readily checked from the relations given in Proposition~\ref{prop:presentation} 
that this defines a $D(B)$-module of dimension $l$. The next result is proved 
in~\cite[Theorem~2.5]{chen 1}.

\bigskip

\begin{theo}[Chen]\label{theo:simples}
The map
$$\fonctio{\{1,2,\dots,d\} \times \ZM/d\ZM}{\Irr(D(B))}{(l,p)}{M_{l,p}}$$
is bijective.
\end{theo}

\bigskip

\subsection{Blocks} 
We put $\L(d)=\ZM/d\ZM \times \ZM/d\ZM$, a set in
canonical bijection with $\{1,2,\dots,d\} \times \ZM/d\ZM$, which 
parametrizes the simple $D(B)$-modules. Given $\l \in \L(d)$, we denote by $M_\l$ 
the corresponding simple $D(B)$-module. We also set $(\ZM/d\ZM)^\#=(\ZM/d\ZM)\setminus\{0\}$ 
and $\L^\#(d)=(\ZM/d\ZM)^\# \times \ZM/d\ZM$. Finally, let $\L^0(d)=\{0\} \times \ZM/d\ZM$ 
be the complement of $\L^\#(d)$ in $\L(d)$. 

Define
$$\fonction{\iota}{\L(d)}{\L(d)}{(l,p)}{(-l,p+l).}$$
We have $\iota^2=\Id_{\L(d)}$ and $\L^0(d)$ is the set of fixed points of $\iota$. 
Given $\LC$ a $\iota$-stable subset of $\L(d)$, we denote by $[\LC/\iota]$ 
a set of representatives of $\iota$-orbits in $\LC$. The next result is proved  
in~\cite[Theorem~2.26]{egst1}.

\bigskip

\begin{theo}[Erdmann-Green-Snashall-Taillefer]\label{theo:blocks}
Let $\l$, $\l' \in \L(d)$. Then 
$M_\l$ and $M_{\l'}$ belong to the same block of $D(B)$ if and only if 
$\l$ and $\l'$ are in the same $\iota$-orbit.
\end{theo}

\bigskip

We have constructed in~\S\ref{sub:theta} a central element, namely $\th$. 
Note that

\equat\label{eq:action-theta}
\text{\it The element $\th$ acts on $M_{l,p}$ by multiplication by $\z^{(p-1)(l+p-1)}$.}
\endequat

\begin{proof}
It is sufficient to compute the action of $\th$ on the vector $e_1^{(l,p)}$. Note that 
$E^i e_1^{(l,p)}=0$ as soon as $i \ge 1$. Therefore, for computing 
$\th e_1^{(l,p)}$ using the formula~(\ref{eq:theta}), 
only the terms corresponding to $i=0$ remain. Consequently,
\eqna
\o_{l,p}(\th)&=&\DS{\frac{1}{d} \sum_{j,k=0}^{d-1} \z^{-jk} \z^{(2p+l-1)k}\z^{(p+l-1)(j-k-1)}}\\
&=& \DS{\frac{\z^{1-l-p}}{d} \sum_{k=0}^{d-1} \z^{pk}\Bigl(\sum_{j=0}^{d-1}\z^{(p+l-1-k)j}\Bigr)}\\
\endeqna
The term inside the big parenthesis is equal to $d$ if $p+l-1-k \equiv 0 \mod d$, and 
is equal to $0$ otherwise. The result follows.
\end{proof}

\bigskip

\subsection{Projective modules} 
Given $\l \in \L(d)$, we denote by $P_\l$ a projective cover of $M_\l$. 
The next result is proved in~\cite[Corollary~2.25]{egst1}.

\bigskip

\begin{theo}[Erdmann-Green-Snashall-Taillefer]\label{theo:projectifs}
Let $\l \in \L(d)$.
\begin{itemize}
\itemth{a} If $\l \in \L^\#(d)$, then $\dim_\CM(P_\l) = 2d$, $\Rad^3(P_\l)=0$ 
and the Loewy structure of $P_\l$ is given by:
$$
\begin{array}{rcc}
P_\l/\Rad(P_\l) &\simeq & M_\l\\
\Rad(P_\l)/\Rad^2(P_\l) &\simeq & M_{\iota(\l)} \oplus M_{\iota(\l)}\\
\Rad^2(P_{\l}) &\simeq& M_{\l}\\
\end{array}
$$

\itemth{b} $P_{d,p}=M_{d,p}$ has dimension $d$.
\end{itemize}
\end{theo}

\bigskip

\section{Tensor structure}

\medskip

We mainly refer here to the work of Erdmann, Green, Snashall and Taillefer~\cite{egst1},~\cite{egst2}. 
Since $D(B)$ is a finite dimensional Hopf algebra, the category 
$D(B)\modules$ inherits a structure of a tensor category. We will 
compute here some tensor products between simple modules. 
We will denote by $M_l$ the simple module $M_{l,0}$.

\bigskip

\subsection{Invertible modules} 
We denote by $V_\xi =\CM v_\xi$ 
the one-dimensional $D(B)$-module associated with the morphism $\e_\xi : D(B) \to \CM$ 
defined in \S\ref{sub:morphismes}: 
$$b v_\xi=\e_\xi(b) v_\xi$$
for all $b \in D(B)$. We have 
\equat\label{eq:v-xi}
V_{\z^p}\simeq M_{1,p}.
\endequat
An immediate computation using the comultiplication $\D$ shows that
\equat\label{eq:tenseur-v-xi}
M_{l,p} \otimes V_{\z^q} \simeq V_{\z^q} \otimes M_{l,p} \simeq M_{l,p+q}
\endequat
as $D(B)$-modules. The $V_\xi$'s are (up to isomorphism) 
the only invertible objects in the tensor category $D(B)\module$.

\bigskip

\subsection{Tensor product with $M_2$}
We set $e_i=e_i^{(2,0)}$ for $i \in \{1,2\}$, so that 
$(e_1,e_2)$ is the standard basis of $M_2$. The next result 
is a particular case of~\cite[Theorem~4.1]{egst1}.

\bigskip

\begin{theo}[Erdmann-Green-Snashall-Taillefer]\label{theo:tenseur-m2}
Let $(l,p)$ be an element of 
$\{1,2,\dots,d\} \times \ZM/d\ZM$.
\begin{itemize}
\itemth{a} If $l \le d-1$, then 
$M_2 \otimes M_{l,p} \simeq M_{l+1,p} \oplus M_{l-1,p+1}$.

\itemth{b} $M_2 \otimes M_{d,p} \simeq P_{d-1,p}$.
\end{itemize}
\end{theo}

\bigskip

\section{Grothendieck rings}

\bigskip

We denote by $\groth(D(B))$ the Grothendieck ring of the 
category of (left) $D(B)$-modules. 
%For simplifying notation, we will change 
% our indexing set for simple $D(B)$-modules: we set 
% $$\L(d)=\ZM/d\ZM \times \ZM/d\ZM,\quad 
% (\ZM/d\ZM)^\#=\ZM/d\ZM \setminus \{0\},$$
% $$\L^0(d)=\{0\} \times \ZM/d\ZM \quad\text{and}\quad
% \L^\#(d)=(\ZM/d\ZM)^\# \times \ZM/d\ZM=\L(d) \setminus \L^0(d).$$
% We denote by $\iota$ the involution of $\L(d)$ defined by 
% $\iota(i,j)=(-i,i+j)$. If $\l=(i,j) \in \L(d)$, we denote by 
% $M_\l$ the simple modules $M_{l,p}$ where $l$ is the unique element 
% of ...

\bigskip

\subsection{Structure}
Since $D(B)$ is a braided 
Hopf algebra (with universal $R$-matrix $R$), 
\equat\label{eq:commutativite}
\text{\it the ring $\groth(D(B))$ is commutative.}
\endequat
Given $M$  a $D(B)$-module, we denote by $\isomorphisme{M}$ the class of $M$ in $\groth(D(B))$. We set
$$\mb_\l=\isomorphisme{M_\l},
\qquad\mb_{l}=\isomorphisme{M_{l,0}}\qquad\text{and}
\qquad\vb_\xi=\isomorphisme{V_\xi} \in \groth(D(B)),$$
% $$\mb_{l,p}^\sta=\isomorphisme{M_{l,p}}_\sta,
% \qquad\mb_{l}^\sta=\isomorphisme{M_{l,0}}_\sta,
% \qquad\vb_\xi^\sta=\isomorphisme{V_\xi}_\sta \in \groth^\sta(D(B)).\leqno{\text{and}}$$
Recall that $\vb_{\z^p}=\mb_{1,p}$.% and $\vb_{\z^p}^\sta=\mb_{1,p}^\sta$. 
% $$\pb_{l,p}=\isomorphisme{P_{l,p}} \in \groth(D(B)),$$
% $$\pb_{l}=\isomorphisme{P_{l,0}} \in \groth(D(B)).\leqno{\text{and}}$$
It follows from~(\ref{eq:tenseur-v-xi}) and Theorem~\ref{theo:tenseur-m2} that
\equat\label{eq:grothendieck}
\vb_{\z^q}\mb_{l,p}=\mb_{l,p+q}\quad\text{and}\quad
\mb_2 \mb_{l,p}=
\begin{cases}
\mb_{l+1,p} + \mb_{l-1,p+1} & \text{if $l \le d-1$,}\\
2(\mb_{d-1,p} + \mb_{1,p-1}) & \text{if $l=d$.}\\
\end{cases}
\endequat

\bigskip

\begin{prop}\label{prop:structure}
The Grothendieck ring $\groth(D(B))$ is generated by $\vb_\z$ and $\mb_2$. 
% The kernel of the morphism $\groth(D(B)) \to \groth^\sta(D(B))$ is 
% the ideal generated by $\mb_d$.
\end{prop}

\bigskip

\begin{proof}
We will prove by induction on $l$ that $\mb_{l,p} \in \ZM[\vb_\z,\mb_2]$. 
Since $\mb_{1,p}=(\vb_\z)^p$, this is true for $l=1$. Since 
$\mb_{2,p}=(\vb_\z)^p \mb_2$, this is also true for $l=2$. Now the induction 
proceeds easily by using~(\ref{eq:grothendieck}). 
%This proves the first assertion.
\end{proof}

\bigskip
\def\ST{{\mathrm{st_B}}}

\def\trace{{\Tb\rb}}

\subsection{Some characters} 
If $b \in D(B)$ is {\it group-like}, then the map
$$\fonctio{\groth(D(B))}{\CM}{\isomorphisme{M}}{\trace(b|_M)}$$
is a morphism of rings. Here, $\trace$ denotes the usual trace (not the quantum trace) 
of an endomorphism of a finite dimensional vector space. 
Recall from Lemma~\ref{lem:group-like} that the only group-like elements of $D(B)$ are 
the $K^i z^j$, where $(i,j) \in \L(d)$. We set
$$\fonction{\chi_{i,j}}{\groth(D(B))}{\CM}{\isomorphisme{M}}{\trace(K^iz^j|_M).}$$
An easy computation yields
\equat\label{eq:chiij}
\chi_{i,j}(\mb_{l,p})=\z^{pi+(2p+l-1)j}\cdot(l)_{\z^i}.
\endequat
Note that the $\chi_{i,j}$'s are not necessarily distinct:

\bigskip

\begin{lem}\label{lem:chiij}
Let $\l$ and $\l'$ be two elements of $\L(d)$. Then 
$\chi_{\l}=\chi_{\l'}$ if and only if $\l$ and $\l'$ are in the same $\iota$-orbit. 
\end{lem}

\bigskip

\begin{proof}
Let us write $\l=(i,j)$ and $\l'=(i',j')$. 
The ``if'' part follows directly from~(\ref{eq:chiij}). Conversely, assume that 
$\chi_{i,j}=\chi_{i',j'}$. By applying these two characters to $\vb_\z$ and $\mb_2$, we get:
$$
\begin{cases}
\z^{i+2j}=\z^{i'+2j'},\\
\z^j(1+\z^i)=\z^{j'}(1+\z^{i'}).
\end{cases}
$$
It means that the pairs $(\z^j,\z^{i+j})$ and $(\z^{j'},\z^{i'+j'})$ have the same sum 
and the same product, so $(\z^j,\z^{i+j})=(\z^{j'},\z^{i'+j'})$ or 
$(\z^j,\z^{i+j})=(\z^{i'+j'},\z^{j'})$. In other words, $(i',j')=(i,j)$ or 
$(i',j')=\iota(i,j)$, as expected.
\end{proof}

\bigskip

\section{Triangulated categories}
\bigskip

\subsection{Stable category} 
As $B$ is a Hopf algebra, it is selfinjective, i.e., $B$ is an injective $B$-module.
Recall that the stable category $B\stable$ of $B$ is 
the additive category
quotient of $B\modules$ by the full subcategory $B\projectifs$ of
projective modules. Since $B$ is selfinjective, the category $B\stable$ has a natural
triangulated structure.
Similarly, the category $D(B)\stable$ is triangulated. Note that 
a $B$-module (resp. a $D(B)$-module) is projective if and only if it is injective. 
Since the tensor product of a projective 
$D(B)$-module by any $D(B)$-module is still projective~\cite[Proposition~4.2.12]{egno}, 
it inherits a structure of monoidal category (such that the 
canonical functor $D(B)\modules \to D(B)\stable$ is monoidal). 
In particular, its Grothendieck group (as a triangulated category), 
which will be denoted by $\groth^{\sta}(D(B))$, is a ring and 
the natural map
$$\fonctio{\groth(D(B))}{\groth^\sta(D(B))}{\mb}{\mb^\sta}$$
is a morphism of rings. Given $M$ a $D(B)$-module, we denote by 
$\isomorphisme{M}_\sta$ its class in $\groth^\sta(D(B))$.

It follows from Theorem~\ref{theo:projectifs} that
\equat\label{eq:grothendieck-stable}
\mb_{d,p}^\sta=0\qquad\text{and}\qquad 2(\mb_{l,p}^\sta + \mb_{d-l,p+l}^\sta)=0
\endequat
if $l \le d-1$. 

\bigskip

\subsection{A further quotient}
We denote by $D(B)\PROJECTIFS$ the full subcategory of $D(B)\modules$ whose objects 
are the $D(B)$-modules $M$ such that $\Res_B^{D(B)} M$ is a projective $B$-module. 
Since $D(B)$ is a free $B$-module (of rank $d^2$), $D(B)\projectifs$ is a full subcategory 
of $D(B)\PROJECTIFS$. 
We denote by $D(B)\STABLE$ the additive 
quotient of the category $D(B)\modules$ by the full 
subcategory $D(B)\PROJECTIFS$: it is also the quotient of $D(B)\stable$ 
by the image of $D(B)\PROJECTIFS$ in $D(B)\stable$. 

\bigskip

\begin{lem}\label{lem:thick}
The image of $D(B)\PROJECTIFS$ in $D(B)\stable$ is a thick triangulated subcategory. 
In particular, $D(B)\STABLE$ is triangulated.
\end{lem}

\bigskip

\begin{proof}
Given $M$ a $D(B)$-module, we denote by $\pi_M : P(M) \surto M$ (resp. $i_M : M \injto I(M)$) 
a projective cover (resp. injective hull) of $M$. We need to prove the following facts:
\begin{itemize}
\itemth{a} If $M$ belongs to $D(B)\PROJECTIFS$, then $\Ker(\pi_M)$ and $I(M)/\im(i_M)$ 
also belong to $D(B)\PROJECTIFS$.

\itemth{b} If $M \oplus N$ belongs to $D(B)\PROJECTIFS$, then $M$ and $N$ also belong to 
$D(B)\PROJECTIFS$.

\itemth{c} If $M$ and $N$ belong to $D(B)\PROJECTIFS$ and $f : M \to N$ is a morphism 
of $D(B)$-modules, then the cone of $f$ also 
belong to $D(B)\PROJECTIFS$.
\end{itemize}

\medskip

(a) Assume that $M$ belongs to $D(B)\PROJECTIFS$. Since it is a projective $B$-module, 
there exists a morphism of $B$-modules $f : M \to P(M)$ such that $\pi_M \circ f = \Id_M$. 
In particular, $P(M) \simeq \Ker(\pi_M) \oplus M$, as a $B$-module. So $\Ker(\pi_M)$ 
is a projective $B$-module. 

On the other hand, $I(M)$ is a projective $D(B)$-module since $D(B)$
so it is a projective $B$-module and so it is an injective $B$-module. 
So, again, $I(M) \simeq M \oplus I(M)/\im(i_M)$, so $I(M)/\im(i_M)$ is a projective 
$B$-module. This proves~(a).

\medskip

(b) is obvious.

\medskip

(c) Let $M$ and $N$ belong to $D(B)\PROJECTIFS$ and $f : M \to N$ be a morphism 
of $D(B)$-modules. Let $\D_f : M \to I(M) \oplus N$, $m \mapsto (i_M(m),f(m))$. Then 
the cone of $f$ is isomorphic in $D(B)\stable$ to $(I(M) \oplus N)/\im(\D_f)$. 
But $\D_f$ is injective, $M$ is an injective $B$-module 
and so $I(M) \simeq M \oplus (I(M) \oplus N)/\im(\D_f)$ as a $B$-module, which shows that 
$(I(M) \oplus N)/\im(\D_f)$ is a projective $B$-module. 
\end{proof}

\bigskip

We denote by $\groth^\ST(D(B))$ the Grothendieck group of
$D(B)\STABLE$, viewed as a triangulated category.
If $M$ belongs to $D(B)\PROJECTIFS$ 
and $N$ is any $D(B)$-module, then $M \otimes N$ and $N \otimes M$ are projective 
$B$-modules~\cite[Proposition~4.2.12]{egno}, so $D(B)\STABLE$ inherits 
a structure of monoidal category, compatible with the triangulated structure. 
This endows $\groth^\ST(D(B))$ with a ring structure. The natural map 
$\groth(D(B)) \to \groth^\ST(D(B))$ will be denoted by $\mb \mapsto \mb^\ST$: 
it is a surjective morphism of rings that
factors through $\groth^\sta(D(B))$. 

\medskip

If $\l \in \L^\#(d)$, then it follows from~\cite[Property~1.4]{egst2} that there exists 
a $D(B)$-module $P_\l^B$ which is projective as a $B$-module and such that there 
is an exact sequence 
$$0 \longto M_{\iota(\l)} \longto P^B_{\l} \longto M_{\l} \longto 0.$$
It then follows that 
\equat\label{eq:mlp}
\mb_\l^\ST + \mb_{\iota(\l)}^\ST=0.
\endequat
Also, we still have 
\equat\label{eq:mdp}
\mb_{d,p}^\ST=0.
\endequat

\medskip

The next theorem follows from~(\ref{eq:mlp}),~(\ref{eq:mdp}) and Proposition~\ref{prop:structure}.

\bigskip

\begin{theo}\label{theo:grothendieck-stable}
The ring $\groth^\ST(D(B))$ is generated by $\vb_\z^\ST$ and $\mb_2^\ST$. 
Moreover, 
$$\groth^\ST(D(B))=\bigoplus_{\l \in [\L^\#(d)/\iota]} \ZM \mb_\l^\ST$$
and $\groth^\ST(D(B))$ is a free $\ZM$-module of rank $d(d-1)/2$.
\end{theo}

\bigskip

Recall that Lemma~\ref{lem:chiij} shows that, through the $\chi_{i,j}$'s, 
only $d(d+1)/2$ different characters of the ring $\groth(D(B))$ have been defined. 
It is not clear if $\CM\groth(D(B))$ is semisimple in general but, for 
$d=2$, it can be checked that it is semisimple (of dimension $4$), 
so that there is a fourth character $\groth(D(B)) \to \CM$ which is 
not obtained through the $\chi_{i,j}$'s. 

Now, a character $\chi : \groth(D(B)) \to \CM$ factors through $\groth^\ST(D(B))$ 
if and only if its kernel contains the $\mb_{\l_0}$'s (where $\l_0$ runs over 
$\L^0(d)$) and the $\mb_{\l}+\mb_{\iota(\l)}$'s (where $\l$ runs over $\L^\#(d)$).  
This implies the following result.

\bigskip

\begin{theo}\label{theo:character-stable}
The character $\chi_\l : \groth(D(B)) \to \CM$ factors through 
$\groth^\ST(D(B))$ if and only if $\l \in \L^\#(d)$. So 
the $(\chi_\l)_{\l \in [\L^\#(d)/\iota]}$ 
are all the characters of $\groth^\ST(D(B))$ and the $\CM$-algebra 
$\CM\groth^\ST(D(B))$ is semisimple.
\end{theo}

\bigskip
\def\forget{{\mathbf{For}}}
\def\can{{\mathrm{can}}}

\subsection{Complements}
Given $\CC$ a monoidal category, we denote by $\Zb(\CC)$ its {\it Drinfeld center} 
(see~\cite[\S{XIII.4}]{kassel}) and we denote by $\forget_\CC : \Zb(\CC) \to \CC$ 
the forgetful functor. 

There is an equivalence between $\Zb(B\modules)$ and
$D(B)\modules$ such that
the forgetful functor becomes the restriction functor $\Res_B^{D(B)}$. 
The canonical functors between these categories will be denoted 
by $\can_\sta^{D(B)} : D(B)\modules \to D(B)\stable$, 
$\can_\sta^B : B\modules \to B\stable$ and 
$\can_\ST : D(B)\modules \to D(B)\STABLE$. 
The functor 
$$\can_\sta^B \circ \Res_B^{D(B)} : D(B)\modules \longto B\stable$$
factors through $\Zb(B\stable)$ (this triangulated category needs to be defined
in a homotopic setting, for example that of stable $\infty$-categories).
We obtain a commutative diagram of functors 
$$\diagram
D(B)\modules \rrto^{\DS{\Res_B^{D(B)}}} \ddto_{\DS{\FC}} && B\modules \ddto^{\DS{\can_\sta^B}} \\
&&\\
\Zb(B\stable) \rrto_{\DS{\forget_{B\stable}}} && B\stable.
\enddiagram$$
Since  any $D(B)$-module that is projective as a $B$-module is sent to the zero object 
of $\Zb(B\stable)$ through $\FC$, the functor $\FC$ factors through $D(B)\PROJECTIFS$ and we obtain
a commutative diagram of functors 
$$\diagram
&D(B)\modules \rrto^{\DS{\Res_B^{D(B)}}} \ddto_{\DS{\FC}} \dlto_{\DS{\can_\ST^{D(B)}}} && 
B\modules \ddto^{\DS{\can_\sta^B}} \\
D(B)\STABLE\drto_{\DS{\overline{\FC}}} &&&\\
&\Zb(B\stable) \rrto_{\DS{\forget_{B\stable}}} && B\stable.
\enddiagram$$

\noindent{\bf Question.} {\it Is $\overline{\FC} : D(B)\STABLE \to \Zb(B\stable)$ 
an equivalence of categories?}

\bigskip

\section{Fusion datum}

\medskip

\subsection{Quantum traces} 
The element $R \in D(B) \otimes D(B)$ defined in~\S\ref{sub:braiding} 
is a universal $R$-matrix which endows $D(B)$ with a structure of braided Hopf algebra. 
The category $D(B)\modules$ is braided as follows: given $M$ and $N$ 
two $D(B)$-modules, the braiding $c_{M,N} : M \otimes N \longiso N \otimes M$ 
is given by
$$c_{M,N}(m \otimes n)= \t(R)(n \otimes m).$$
Recall that $\t : D(B) \otimes D(B) \longiso D(B) \otimes D(B)$ is given 
by $\t(a \otimes b)=b \otimes a$. In particular, 
\equat\label{eq:cmn-cnm}
\text{\it $c_{N,M} c_{M,N}: M \otimes N \longiso M \otimes N$ 
is given by the action of $\t(R)R$.}
\endequat

Given $i \in \ZM$, we have $S^2(b)=(z^{-i}K)b(z^{-i}K)^{-1}$ for all $b \in D(B)$ and 
$z^{-i}K$ is group-like, so 
the algebra $D(B)$ is {\it pivotal} with pivot $z^{-i}K$. This endows 
the tensor category $D(B)\modules$ with a structure of pivotal 
category (see Appendix~\ref{sec:preliminaires}) whose associated 
traces $\Tr_+^{(i)}$ and $\Tr_-^{(i)}$ are given as follows: given $M$ a $D(B)$-module 
and $f \in \End_{D(B)}(M)$, we have
$$\Tr_+^{(i)}(f)=\trace(z^{-i}K f)\qquad\text{and}\qquad \Tr_-^{(i)}(f)=\trace(f K^{-1}z^i).$$
Recall that $\trace$ denotes the ``classical'' trace for endomorphisms of 
a finite dimensional vector space. So the pivotal structure depends on the 
choice of $i$ (modulo $d$). The corresponding twist is $\th_i=z^{i}\th$, 
which endows $D(B)\modules$ with a structure of balanced braided category 
(depending on $i$).

\bigskip

\boitegrise{\noindent{\bf Hypothesis and notation.} {\it From now on, 
and until the end of this paper, we assume that 
the Hopf algebra $D(B)$ is endowed with the pivotal structure whose pivot 
is $z^{-1}K$. The structure of balanced braided category is given 
by $\th_{1}=z\th$ and the associated quantum traces $\Tr_\pm^{(1)}$ are denoted 
by $\Tr_\pm$.\\
\hphantom{AA} Given $M$ a $D(B)$-module, we set
$\dim_\pm(M)=\Tr_\pm(\Id_M)$.}}{0.75\textwidth}

\bigskip

We define
$$\dim(D(B))=\sum_{M \in \Irr D(B)} \dim_-(M)\dim_+(M).$$
We have
\equat\label{eq:dim-db-mod}
\dim(D(B))=\frac{2d^2}{(1-\z)(1-\z^{-1})}.
\endequat
This follows easily from the fact that 
\equat\label{eq:dim-mlp}
\dim_+ M_{l,p} = \z^{1-l-p}(l)_\z\quad\text{and}\quad
\dim_- M_{l,p}=\z^{p+l-1}(l)_{\z^{-1}}=\z^p(l)_\z.
\endequat

\bigskip

\subsection{Characters of $\groth(D(B))$ via the pivotal structure} 
As in Appendix~\ref{sec:preliminaires}, these structures (braiding, pivot) 
allow to define characters of $\groth(D(B))$ associated with simple modules (or bricks). 
Given $\l \in \L(d)$,  we set 
$$\fonction{s_{M_\l}^+}{\groth(D(B))}{\CM}{\isomorphisme{M}}{(\Id_{M_{\l}} 
\otimes \Tr_+^M)(c_{M,M_{\l}} c_{M_{\l},M})}$$
$$\fonction{s_{M_{\l}}^-}{\groth(D(B))}{\CM}{\isomorphisme{M}}{(\Tr_-^M
\otimes \Id_{M_{\l}})(c_{M,M_{\l}} c_{M_{\l},M})}\leqno{\text{and}}$$
These are morphism of rings (see Proposition~\ref{prop:caractere}). 
The main result of this section is the following.

\bigskip

\begin{theo}\label{theo:s-matrice}
Given $\l \in \L(d)$, we have
$$s_{M_{\l}}^+=\chi_{-\l}\qquad\text{and}\qquad
s_{M_{(l,p)}}^-=\chi_{(0,1)-\l}.$$
\end{theo}

\bigskip

\def\coef{{\mathrm{coef}}}

\begin{proof}
Write $\l=(l,p)$ and
$$\gamb_{i,j,k}=\frac{\z^{(i-k)(i+j)-i(i+1)/2}}{(i)!_\z}.$$
We have
$$\t(R)R=\frac{1}{d^2}\sum_{i,i',j,j',k,k'=0}^{d-1} \gamb_{i,j,k}\gamb_{i',j',k'} 
\z^{i(k'-j')} (z^{-k'}F^{i'}E^iK^{k'+j}) \otimes (z^{-k}E^{i'}F^iK^{j'+k}).$$
We need to compute the endomorphism of $M_{l,p}$ given by
$$(\Id_{M_{l,p}} \otimes \Tr_+^M)(\t(R)R|_{M_{l,p} \otimes M}).$$
Since $M_{l,p}$ is simple, this endomorphism is the multiplication 
by a scalar $\varpi$, and so it is sufficient to compute the action 
on $e_1^{(l,p)} \in M_{l,p}$. Therefore, all the terms (in the big sum giving 
$\t(R)R$) corresponding to $i \neq 0$ disappear (because $Ee_1^{(l,p)}=0$). 
Also, since we are only interested in the coefficient on $e_1^{(l,p)}$ 
of the result (because the coefficients on other vectors will be zero), 
all the terms corresponding to $i' \neq 0$ also disappear. Therefore,
$$\varpi=\frac{1}{d^2}\sum_{j,j',k,k' \in \ZM/d\ZM} \z^{-kj-k'j'}(\z^{2p+l-1})^{-k'}(\z^{p+l-1})^{k'+j} 
\trace(z^{-1}Kz^{-k}K^{j'+k}|_M).$$
So it remains to compute the element
$$b=\frac{1}{d^2}\sum_{j,j',k,k' \in \ZM/d\ZM} \z^{-kj-k'j'}(\z^{2p+l-1})^{-k'}(\z^{p+l-1})^{k'+j} 
z^{-k-1}K^{j'+k+1}$$
of $D(B)$. Since
$$b=\frac{1}{d^2}\sum_{j',k \in \ZM/d\ZM}\Bigl(\sum_{j,k' \in \ZM/d\ZM} 
\z^{l(p+l-1-k)}\z^{k'(-p-j')}\Bigr) z^{-k-1}K^{j'+k+1},$$
only the terms corresponding to $k=l+p-1$ and $j'=-p$ remain, hence
$$b=z^{-l-p}K^l.$$
So $s_{M_{\l}}^+=\chi_{-\iota(\l)}=\chi_{-\l}$, as expected.

The other formula is obtained via a similar computation.
\end{proof}

\bigskip

% Again, let $\EC$ be a set of representatives of the $\iota$-orbits in 
% $\{1,2,\dots,d-1\} \times \ZM/d\ZM$. 
%We first start by defining a kind of $S$-matrix associated with $D(B)\modules$. 
We denote by 
$\SM^\pm=(\SM_{\l,\l'}^\pm)_{\l,\l' \in \L(d)}$ the 
square matrix defined by
$$\SM_{\l,\l'}^\pm = \Tr_\pm(c_{M_{\l'},M_{\l}} \circ c_{M_{\l},M_{\l'}}).$$
Similary, we define $\TM^\pm$ to be the diagonal matrix (whose rows and columns 
are indexed by $\l \in \L(d)$) and whose $\l$-entry is
$$\TM_{\l}^\pm = \o_{\l}(\th_1^{\mp 1}).$$
Let us first give a formula for $\SM_{\l,\l'}^\pm$ and $\TM_{\l}^\pm$.

\bigskip

\begin{coro}\label{coro:fusion-datum-db}
Let $(l,p)$, $(l',p') \in \L(d)$. We have
$$\SM_{(l,p),(l',p')}^+=\frac{\z}{1-\z} \z^{-ll'-lp'-pl'-2pp'}(1-\z^{ll'}),\qquad
\TM_{(l,p)}^+=\z^{-p(p+l)},$$
$$\SM_{(l,p),(l',p')}^-=\frac{\z^{2p+l+2p'+l'-1}}{1-\z}\z^{-ll'-lp'-pl'-2pp'}(1-\z^{ll'})
\qquad\text{and}\qquad
\TM_{(l,p)}^-=\z^{p(p+l)}.$$
\end{coro}

\bigskip

\begin{proof}
This follows immediately from formulas~(\ref{eq:action-theta}),~(\ref{eq:chiij}),~(\ref{eq:dim-mlp}) 
and Theorem~\ref{theo:s-matrice}.
\end{proof}

% \begin{rema}\label{rem:chiij}
% If $(l,p)$ and $(l',p')$ belong to $\{1,2,\dots,d\} \times \ZM/d\ZM$, then
% \equat\label{eq:chiij-bis}
% \chi_{l,p}(\mb_{l',p'}) = \z^{lp'+p(2p'+l'-1)}(l')_{\z^l}.
% \endequat
% In particular, $\chi_{l,p}=\chi_{-l,p+l}$. This shows that
% $$s_{M_{l,p}}^+=s_{M_{d-l,l+p}}^+,$$
% as expected from ???.
% 
% Moreover, $\chi_{l,p}=\chi_{l',p'}$ if and only if ... Indeed, 
% ...\finl
% \end{rema}
% 
% \bigskip
% 
% \begin{coro}\label{coro:chiij-stable}
% The morphism of rings $s_{M_{l,p}}^+$ (respectively $s_{M_{l,p}}^-$) 
% factorizes through $\groth^\sta_*(D(B))$ if and only if 
% $l < d$.
% \end{coro}
% 
% \bigskip
% 
% Therefore, we have constructed $d(d-1)/2$ different morphisms of rings 
% (because $s_{M_{l,p}}^+=s_{M_{d-l,l+p}}^+$ by Remark~\ref{rem:chiij}) 
% $\groth^\sta_*(D(B)) \longto \CM$. Since the $\ZM$-rank of $\groth_*^\sta(D(B))$ 
% is $d(d-1)/2$, it means that we have obtained all such morphisms:
% 
% \bigskip
% 
% \begin{theo}\label{theo:s-matrice-stable}
% The algebra $\CM\groth_*^\sta(D(B))$ is commutative, semisimple and 
% any character of $\CM\groth_*^\sta(D(B))$ is the $\CM$-linear extension 
% of some $s_{M_{l,p}}^+$ (or some $s_{M_{l,p}}^-$).
% \end{theo}
% 
% \bigskip

\bigskip

\subsection{Fusion datum associated with $D(B)\STABLE$}\label{sub:fusion-db}
Let $\EC$ denote a set of representatives of $\iota$-orbits in $\{1,2,\dots,d-1\} \times \ZM/d\ZM$. 
We define 
$$\dim^\ST(D(B))=\sum_{(l,p) \in \EC} \dim_-(M_{l,p})\dim_+(M_{l,p}).$$
This can be understood in terms of super-categories, as explained
recently by Lacabanne~\cite{abel 1}. We have
$$\dim^\ST(D(B))=\frac{1}{2}\dim(D(B))=\frac{d^2}{(1-\z)(1-\z^{-1})}.$$
So $\dim^\ST D(B)$ is a positive real number and we denote by $\sqrt{\dim^\ST D(B)}$ its 
positive square root. Since $1-\z^{-1}=-\z^{-1}(1-\z)$, there exists a unique 
square root $\sqrt{-\z}$ of $-\z$ such that 
$$\sqrt{\dim^\ST(D(B))} = \frac{d\sqrt{-\z}}{1-\z}.$$
We denote by $\SM^\ST=(\SM_{\l,\l'})_{\l,\l' \in \EC}$ the square matrix 
defined by
$$\SM_{\l,\l'}^\ST=\frac{\SM^+_{\l,\l'}}{\sqrt{\dim^\ST(D(B))}}.$$
We denote by $\TM^\ST$ the diagonal matrix whose $\l$-entry is $\TM_\l^+$ (for $\l \in \EC$). 
It follows from Corollary~\ref{coro:fusion-datum-db} that 
\equat\label{eq:s-matrice-stable}
\SM_{(l,p),(l',p')}^\ST=\frac{\sqrt{-\z}}{d} \z^{-ll'-lp'-pl'-2pp'}(\z^{ll'}-1)\quad\text{and}\quad 
\TM_{(l,p)}^\ST=\z^{-p(l+p)}.
\endequat
The root of unity $\sqrt{-\z}$ appearing in this formula has been interpreted 
in terms of super-categories by Lacabanne~\cite{abel 1}: it is due to the 
fact that our category is not spherical. Finally, note that
\equat\label{eq:--matrice-iota}
\SM_{(l,p),(l',p')}^\ST=-\SM_{\iota(l,p),(l',p')}^\ST=
-\SM_{(l,p),\iota(l',p')}^\ST=\SM_{\iota(l,p),\iota(l',p')}^\ST.
\endequat

\bigskip

\section{Comparison with Malle $\ZM$-fusion datum}

\medskip

We refer to~\cite{malle} and~\cite{cuntz} for most of the material 
of this section. We denote by $\EC(d)$ the set of pairs $(i,j)$ of integers with
$0 \le i < j \le d-1$.

\makeatletter
\newcommand{\extp}{\@ifnextchar^\@extp{\@extp^{\,}}}
\def\@extp^#1{\mathop{\bigwedge\nolimits^{\!\!#1}}}
\makeatother

\subsection{Set-up}
Let $Y=\{0,1,\dots,d\}$ and let $\pi : Y \to \{0,1\}$ be the map defined by
$$\pi(i)=
\begin{cases}
1 & \text{if $i \in \{0,1\}$},\\
0 & \text{if $i \ge 2$.}
\end{cases}
$$
We denote by $\Psi(Y,\pi)$ the set of maps $f : Y \to \{0,1,\dots,d-1\}$ 
such that $f$ is strictly increasing on $\pi^{-1}(0)=\{2,3,\dots,d\}$ 
and strictly increasing on $\pi^{-1}(1)=\{0,1\}$. Since $f$ is injective 
on $\{2,3,\dots,d\}$, there exists a unique element $\kb(f) \in \{0,1,\dots,d-1\}$ 
which does not belong to $f(\{2,3,\dots,d\})$. Note that, since $f$ is strictly increasing 
on $\{2,3,\dots,d\}$, the element
$\kb(f)$ determines the restriction of $f$ to $\{2,3,\dots,d\}$. 
So the map
\equat\label{eq:bij-psi}
\fonctio{\Psi(Y,\pi)}{\EC(d) \times \{0,1,\dots,d-1\}}{f}{(f(0),f(1),\kb(f))}
\endequat
is bijective. For $f \in \Psi(Y,\pi)$, we set 
$$\e(f)=(-1)^{|\{(y,y') \in Y \times Y~|~y<y'~\text{and}~f(y)<f(y')\}|}.$$

\medskip
\def\frob{{\mathrm{Fr}}}

We put by $V=\bigoplus_{i=0}^{d-1} \CM v_i$
and we denote by $\SC$ the square matrix $(\z^{ij})_{0 \le i,j \le d-1}$, 
which will be viewed as an automorphism of $V$. Note that $\SC$ is the character stable 
of the cyclic group $\mub_d$. 
We set $\d(d)= \det(\SC)=\prod_{0 \le i < j \le d-1} (\z^j-\z^i)$. Recall that
$$\d(d)^2=(-1)^{(d-1)(d-2)/2} d^d.$$
Given $f \in \Psi(Y,\pi)$, let 
$$\vb_f = (v_{f(0)} \wedge v_{f(1)}) \otimes (v_{f(2)} \wedge v_{f(3)} \wedge \cdots \wedge 
v_{f(d)}) \in \bigl(\extp^2 V\bigr) \otimes \bigl(\extp^{d-1} V\bigr).$$
Note that $(\vb_f)_{f \in \Psi(Y,\pi)}$ is a $\CM$-basis of 
$\bigl(\extp^2 V\bigr) \otimes \bigl(\extp^{d-1} V\bigr)$. Given $f' \in \Psi(Y,\pi)$, 
we put
$$\Bigl(\bigl(\extp^2 \SC\bigr) \otimes \bigl(\extp^{d-1} \SC\bigr)\Bigr)(\vb_{f'})=
\sum_{f \in \Psi(Y,\pi)} \Sb_{f,f'} \vb_f.
$$
In other words, $(\Sb_{f,f'})_{f,f' \in \Psi(Y,\pi)}$ is the matrix of the automorphism 
$\bigl(\extp^2 \SC\bigr) \otimes \bigl(\extp^{d-1} \SC\bigr)$ of 
$\bigl(\extp^2 V\bigr) \otimes \bigl(\extp^{d-1} V\bigr)$ in the basis 
$(\vb_f)_{f \in \Psi(Y,\pi)}$.

\bigskip

\begin{lem}\label{lem:formule-s}
Let $f$, $f' \in \Psi(Y,\pi)$. We define
$$i=f(0),\quad j=f(1),\quad k=\kb(f),$$
$$i'=f'(0),\quad j'=f'(1),\quad k'=\kb(f').$$
We have
$$\Sb_{f,f'}=(-1)^{k+k'}\frac{\d(d)}{d} \z^{-kk'}(\z^{ii'+jj'}-\z^{ij'+ji'}).$$
\end{lem}

\bigskip

\begin{proof}
The computation of the action of $\extp^2 \SC$ is easy, and gives the term 
$\z^{ii'+jj'}-\z^{ij'+ji'}$. It remains to show that the determinant 
of the matrix $\SC(k,k')$ obtained from $\SC$ by removing the $k$-th row and the $k'$-th column 
is equal to $(-1)^{k+k'}\z^{-kk'} \d(d)/d$. For this, let $\SC'(k)$ denote 
the matrix whose $k$-th row is equal to $(1,t,t^2,\dots,t^{d-1})$ (where $t$ is an 
indeterminate) and whose other rows coincide with those of $\SC$. Then 
$(-1)^{k+k'}\det(\SC(k,k'))$ is equal to the coefficient of $t^{k'}$ in the polynomial 
$\det(\SC'(k))$. This is a Vandermonde determinant and
\eqna
\det(\SC'(k))&=&\DS{\prod_{\substack{0 \le i < j \le d-1 \\ i \neq k, j \neq k}}(\z^j-\z^i) \cdot 
\prod_{i =0}^{k-1}(t-\z^i) \cdot \prod_{i=k+1}^{d-1}(\z^i-t)}\\
&=& \d(d) \prod_{\substack{i =0 \\ i \neq k}}^{d-1} \frac{(t-\z^i)}{(\z^k-\z^i)}.
\endeqna
Since
$$\prod_{\substack{i =0 \\ i \neq k}}^{d-1} (t-\z^i)=\frac{t^d-1}{t-\z^k}= 
\sum_{i=0}^{d-1} t^i \z^{(d-1-i)k},$$
we have
$$\det(\SC(k,k'))=(-1)^{k+k'}\d(d)\frac{\z^{(d-1-k')k}}{d\z^{(d-1)k}}=
(-1)^{k+k'}\frac{\d(d)}{d} \z^{-kk'},$$
as desired.
\end{proof}

\bigskip

\subsection{Malle $\ZM$-fusion datum} 
Let
$$\Psi^\#(Y,\pi)=\{f \in \Psi(Y,\pi)~|~\sum_{y \in Y} f(y) \equiv \frac{d(d-1)}{2} \mod d\}.$$
Given $f \in \Psi^\#(Y,\pi)$, we define
$$\frob(f)=\z_*^{d(1-d^2)}\prod_{y \in Y} \z_*^{-6(f(y)^2+df(y))},$$
where $\z_*$ is a primitive $(12d)$-th root of unity such that $\z_*^{12}=\z$.

\bigskip

We denote by $\TM$ diagonal matrix (whose rows and column are indexed by $\Psi^\#(Y,\pi)$) 
equal to $\diag(\frob(f))_{f \in \Psi^\#(Y,\pi)}$. We denote by 
$\SM=(\SM_{f,g})_{f,g \in \Psi^\#(Y,\pi)}$ the square matrix 
defined by
$$\SM_{f,g}=\frac{(-1)^{d-1}}{\d(d)} \e(f)\e(g)~\overline{\Sb}_{f,g}.$$
Note that $\SM_{f,f_{0,1}} \neq 0$ for all $f \in \Psi^\#(Y,\pi)$ (see Lemma~\ref{lem:formule-s}). 

\bigskip

\begin{prop}[Malle~\cite{malle}, Cuntz~\cite{cuntz}]\label{prop:fusion-datum}
With the previous notation, we have:
\begin{itemize}
\itemth{a} $\SM^4=(\SM\TM)^3=[\SM^2,\TM]=1$.

\itemth{b} $\lexp{t}{\SM}=\SM$ and $\lexp{t}{\overline{\SM}}~ \SM=1$.

\itemth{c} For all $f$, $g$, $h \in \Psi^\#(Y,\pi)$, the number
$$N_{f,g}^h = 
\sum_{i \in \Psi^\#(Y,\pi)} \frac{\SM_{i,f}\SM_{i,g}\overline{\SM}_{i,h}}{\SM_{i,f_{0,1}}}$$
belongs to $\ZM$. 
\end{itemize}
\end{prop}

\bigskip

The pair $(\SM,\TM)$ is called the {\it Malle $\ZM$-fusion datum}. 

\bigskip

\subsection{Comparison}
We wish to compare the $\ZM$-fusion datum $(\SM,\TM)$ with the ones 
obtained from the tensor categories $D(B)\modules$ and $D(B)\stable$. For this, we will 
use the bijection~(\ref{eq:bij-psi}) to characterize elements of $\Psi^\#(Y,\pi)$.
Given $k \in \ZM$, we denote by $k^\res$ the unique element in $\{0,1,\dots,d-1\}$ 
such that $k \equiv k^\res \mod d$.

\bigskip

\begin{lem}\label{lem:fij}
Let $f \in \Psi(Y,\pi)$. Then $f \in \Psi^\#(Y,\pi)$ if and only if $\kb(f)=(f(0)+f(1))^\res$.
Consequently, the map
$$\fonctio{\Psi^\#(Y,\pi)}{\EC(d)}{f}{(f(0),f(1))}$$
is bijective
\end{lem}

\bigskip

\begin{proof}
We have
$$\sum_{y \in Y} f(y) = f(0)+f(1) + \frac{d(d-1)}{2}-\kb(f)$$
and the result follows.
\end{proof}

\bigskip

Given $(i,j) \in \EC(d)$, we denote by $f_{i,j}$ the unique element of $\Psi^\#(Y,\pi)$ 
such that $f_{i,j}(0)=i$ and $f_{i,j}(1)=j$. We have
\equat\label{eq:frob-fij}
\frob(f_{i,j})=\z^{ij}
\endequat
and, if $(i,j)$, $(i',j') \in \L(d)$, then
\equat\label{eq:s-fij}
\SM_{f_{i,j},f_{i',j'}}=\frac{(-1)^{(i+j)^\res+(i'+j')^\res}}{d}\e(f_{i,j})\e(f_{i',j'})
(\z^{ij'+ji'}-\z^{ii'+jj'}).
\endequat

\begin{proof}
The second equality follows immediately from Lemmas~\ref{lem:formule-s} and~\ref{lem:fij}.
Let us prove the first one. 
By definition, $\frob(f_{i,j})=\z_*^\a$, where 
$$\a=d(1-d^2)-6 \sum_{y \in Y} (f_{i,j}(y)^2 + d f_{i,j}(y)).$$
The construction of $f_{i,j}$ shows that 
$$\a=d(1-d^2)-6(i^2+di)-6(j^2+dj)-6 \sum_{k=0}^{d-1}(k^2+dk) + 6(((i+j)^\res)^2+d(i+j)^\res).$$
Write $i+j=(i+j)^\res + \eta d$, with $\eta \in \{0,1\}$. Then $\eta^2=\eta$ and so 
\eqna
(i+j)^2+ d(i+j) &=&((i+j)^\res)^2+d(i+j)^\res + 2\eta d(i+j) + 2\eta d^2 \\
&\equiv& 
((i+j)^\res)^2+d(i+j)^\res \mod 2d.
\endeqna
Therefore, 
\eqna
\a &\equiv& \DS{12ij +d(1-d^2)  -6 \sum_{k=0}^{d-1}(k^2+dk) \mod 12d}\\
&\equiv& 12ij \mod 12d.
\endeqna
So $\frob(f_{i,j})=\z_*^{12ij}=\z^{ij}$.
\end{proof}

\bigskip

We define
$$\fonction{\ph}{\EC(d)}{\L^\#(d)}{(i,j)}{(j-i,i).}$$
Note that $\ph(\EC(d))$ is a set of representatives of $\iota$-orbits in 
$\L^\#(d)$. We set
$$\pht(i,j)=
\begin{cases}
\ph(i,j) & \text{if $(-1)^{(i+j)^\res}\e(f_{i,j})=1$,} \\
\iota(\ph(i,j)) & \text{if $(-1)^{(i+j)^\res}\e(f_{i,j})=-1$.}
\end{cases}
$$
Then $\pht(\EC(d))$ is also a set of representatives of $\iota$-orbits in 
$\L^\#(d)$ and the pairs of matrices $(\SM^\ST,\TM^\ST)$ and $(\SM,\TM)$ 
are related by the following equality (which follows immediately from 
Corollary~\ref{coro:fusion-datum-db} and formulas~(\ref{eq:--matrice-iota}),~(\ref{eq:frob-fij}) 
and~(\ref{eq:s-fij})):
\equat\label{eq:comparaison}
\overline{\SM}_{f_{i,j},f_{i',j'}}=\sqrt{-\z} ~\SM^\ST_{\pht(i,j),\pht(i',j')}\qquad\text{and}\qquad
\overline{\TM}_{f_{i,j}}=\TM^\ST_{\pht(i,j)}.
\endequat
Therefore, up to the change of $\z$ into $\z^{-1}$, we obtain our main result.

\bigskip

\begin{theo}\label{theo:main}
Malle $\ZM$-fusion datum $(\SM,\TM)$ can be categorified by the monoidal 
category $D(B)\STABLE$, endowed with the pivotal structure induced by the pivot 
$z^{-1}K$ and the balanced structure induced by $z\th$.
\end{theo}

%\bigskip
%
%\begin{rema}\label{rem:exceptions}
%M. Brou\'e, G. Malle \& J. Michel have associated to a class of exceptional 
%reflection groups (the {\it spetsial} ones) a set of {\it ``unipotent 
%characters''} and a partition of these unipotent characters into families. 
%To each family, they have also associated a $\ZM$-fusion datum $(S,T)$. 
%It turns out that some of these $\ZM$-fusion data can by categorified 
%by Hopf quotients of the algebra $D(B)$ (see~\cite{abel 2}).\finl
%\end{rema}

\setcounter{section}{0}
\renewcommand\thesection{\Alph{section}}

\section{Appendix. Reminders on $S$-matrices}\label{sec:preliminaires}

\medskip

We follow closely~\cite[Chapters~4~and~8]{egno}.

Let $\CC$ be a {\it tensor category} over $\CM$, as defined in~\cite[Definition~4.1.1]{egno}: 
$\CC$ is a locally finite $\CM$-linear rigid monoidal category (whose unit object 
is denoted by $\unb$) such that the bifunctor $\otimes : \CC \times \CC \to \CC$ 
is $\CM$-bilinear on morphisms and $\End_\CC(\unb)=\CM$. 
If $X$ is an object in $\CC$, its left (respectively right) dual is denoted 
by $X^*$ (respectively $\lexp{*}{X}$) and we denote by 
$$\coev_X : \unb \longto X \otimes X^*\qquad\text{and}\qquad
\ev_X : X^* \otimes X \longto \unb$$
the {\it coevealuation} and {\it evaluation} morphisms respectively. 

We assume that $\CC$ is {\it braided}, namely that it is endowed 
with a bifunctorial family of isomorphisms $c_{X,Y} : X \otimes Y \longiso Y \otimes X$ 
such that 
\equat\label{eq:cxyy}
c_{X, Y \otimes Y'} = (\Id_Y \otimes c_{X,Y'}) \circ (c_{X,Y} \otimes \Id_{Y'})
\endequat
and
\equat\label{eq:cxxy}
c_{X \otimes X',Y} = (c_{X,Y} \otimes \Id_{X'}) \circ (\Id_X \otimes c_{X',Y}).
\endequat
for all objects $X$, $X'$, $Y$ and $Y'$ in $\CC$ (we have omitted 
the associativity constraints).

Finally, we also assume that $\CC$ is {\it pivotal}~\cite[Definition~4.7.8]{egno}, 
i.e. that it is equipped with a family 
of functorial isomorphisms $a_X : X \to X^{**}$ (for $X$ running over 
the objects of $\CC$) such that $a_{X \otimes Y} = a_X \otimes a_Y$. 
Given $f \in \End_\CC(X)$, the pivotal structure allows to define two {\it traces}:
$$\Tr_+(f) = \ev_{X^*} \circ (a_Xf \otimes \Id_{X^*}) \circ  \coev_X
\in \End_\CC(\unb)=\CM$$
$$\Tr_-(f) = \ev_{X} \circ (\Id_{X^*} \otimes f a_X^{-1}) \circ  \coev_{X^*}
\in \End_\CC(\unb)=\CM.\leqno{\text{and}}$$
We will sometimes write $\Tr_+^X(f)$ or $\Tr_-^X(f)$ for $\Tr_+(f)$ and $\Tr_-(f)$. 
We define two {\it dimensions} 
$$\dim_+(X)=\Tr_+(\Id_X)\qquad\text{and}\qquad \dim_-(X)=\Tr_-(\Id_X).$$
To summarize, we will work under the following hypothesis:

\medskip

\boitegrise{{\bf Hypothesis and notation.} {\it We fix in this section a 
braided pivotal tensor category $\CC$ as above. We denote by $\groth(\CC)$ 
its Grothendieck ring. Given $X$ is an object in $\CC$, we denote by 
$\isomorphisme{X}$ its class in $\groth(\CC)$. The set of isomorphism 
classes of simple objects in $\CC$ will be denoted by $\Irr(\CC)$. 
If $X \in \Irr(\CC)$ and $Y$ is any object in $\CC$, we denote by 
$[Y : X]$ the multiplicity of $X$ in a Jordan-H\"older series of $Y$.}}{0.75\textwidth}

\bigskip

%\subsection{$S$-matrices} 
Given $X$, $Y$ two objects in $\CC$, we set 
$$s_{X,Y}^+ = (\Id_X \otimes \Tr_+^Y)(c_{Y,X} c_{X,Y}) \in \End_\CC(X).$$
$$s_{X,Y}^- = (\Tr_-^Y \otimes \Id_X)(c_{Y,X} c_{X,Y}) \in \End_\CC(X).
\leqno{\text{and}}$$
These induce two morphisms of abelian groups
$$\fonction{s_X^+}{\groth(\CC)}{\End_\CC(X)}{\isomorphisme{Y}}{s_{X,Y}^+}
\qquad\text{and}\qquad
\fonction{s_X^-}{\groth(\CC)}{\End_\CC(X)}{\isomorphisme{Y}}{s_{X,Y}^-.}$$

\bigskip

\begin{defi}\label{defi:endominimal}
An object $X$ in $\CC$ is called a {\bfit brick} if $\End_\CC(X)=\CM$.
\end{defi}

\bigskip

For instance, a simple object is a brick (and $\unb$ is also a brick, 
but $\unb$ is simple in a tensor category~\cite[Theorem~4.3.1]{egno}). 
Note also that a brick is indecomposable. 
So if $\CC$ is moreover semisimple, then an object is a brick if and only 
if it is simple.

If $X$ is a brick, then we will view $s_{X,Y}^+$ and $s_{X,Y}^-$ as elements 
of $\CM=\End_\CC(X)$. 

\bigskip

\begin{prop}\label{prop:caractere}
If $X$ is a brick, then $s_X^+ : \groth(\CC) \to \CM$ 
and $s_X^- : \groth(\CC) \to \CM$ are morphisms of rings.
\end{prop}

\bigskip

\begin{proof}
Assume that $X$ is a brick. We will only prove the result for 
$s_X^+$, which amounts to show that 
$$s_{X,Y \otimes Y'}^+ = s_{X,Y}^+ s_{X,Y'}^+.\leqno{(*)}$$
First, note that the following equality 
$$c_{Y \otimes Y',X} c_{X,Y \otimes Y'} = 
(c_{Y,X} \otimes \Id_{Y'}) \circ (\Id_Y \otimes c_{Y',X}c_{X,Y'}) 
\circ (c_{X,Y} \otimes \Id_{Y'})$$
holds by~(\ref{eq:cxyy}) and~(\ref{eq:cxxy}). Taking 
$\Id_X \otimes \Id_Y \otimes \Tr_+^{Y'}$ on the right-hand side, one gets
$s_{X,Y'}^+ c_{Y,X}c_{X,Y} \in \End_\CC(X \otimes Y)$ (because $X$ is a brick). 
Applying now $\Id_X \otimes \Tr_+^Y$, one get $s_{X,Y'}^+s_{X,Y}^+ \Id_X$. 
Since 
$$(\Id_X \otimes \Tr_+^Y) \circ (\Id_X \otimes \Id_Y \otimes \Tr_+^{Y'})
=\Id_X \otimes \Tr_+^{Y\otimes Y'},$$
this proves $(*)$.
\end{proof}

\bigskip

% \begin{prop}\label{prop:dim-s}
% If $X$ and $Y$ are endominimal, then 
% $$\dim_-(X) s_{X,Y}^+ = \dim_+(Y) s_{X,Y}^-.$$
% \end{prop}
% 
% \bigskip
% 
% \begin{proof}
% This follows immediately from the fact that 
% $$\Tr_-^X \otimes \Tr_+^Y = \Tr_-^X \circ (\Id_X \otimes \Tr_+^Y) = 
% \Tr_+^Y \circ (\Tr_-^X \otimes \Id_Y),$$
% by applying these equalities to $c_{Y,X}c_{X,Y}$.
% \end{proof}
% 
% \bigskip

\begin{prop}\label{prop:subquotient}
Let $X$ be a brick and let $X'$ be a subquotient of $X$ which is also a brick. Then 
$$s_X^+=s_{X'}^+\qquad\text{and}\qquad s_X^-=s_{X'}^-.$$
\end{prop}

\bigskip

\begin{proof}
Indeed, the endomorphism $(\Id_X \otimes \Tr_+^Y)(c_{Y,X} c_{X,Y})$ 
of $X$ is multiplication by a scalar, and this scalar can be computed 
on any non-trivial subquotient of $X$.
\end{proof}

\bigskip

\begin{coro}\label{coro:blocks}
Let $X$ and $X'$ be two bricks in $\CC$ belonging to the same 
block. Then 
$$s_X^+=s_{X'}^+\qquad\text{and}\qquad s_X^-=s_{X'}^-.$$
\end{coro}

\bigskip

\begin{proof}
By Proposition~\ref{prop:subquotient}, we may asume that $X$ and $X'$ are simple. 
We may also assume that $X$ is not isomorphic to $X'$ and that 
$\Ext^1_\CC(X,X')=0$. Let $\Xb \in \CC$ such that 
there exists a non-split exact sequence 
$$0 \longto X' \longto \Xb \longto X \longto 0.$$
Since $X \not\simeq X'$, we have $\End_\CC(\Xb)=\CM$, hence $\Xb$ is a brick. 
It follows from Proposition~\ref{prop:subquotient} that
$$s_\Xb^+=s_X^+=s_{X'}^+\qquad\text{and}\qquad s_\Xb^-=s_X^-=s_{X'}^-,$$
as desired.
\end{proof}

\end{document}